\documentclass{commat}

\title{Asymptotic formula for the multiplicative function $\frac{d(n)}{k^{\omega(n)}}$}

\author{Meselem Karras}

\affiliation{
    \address{Faculty of Science and Technology Djilali Bounaama Khemis Miliana
    University, 44225, Algeria}
    \email{karras.m@hotmail.fr}
    }

\abstract{
For a~fixed integer $k$, we define the multiplicative function
\[
D_{k,\omega}(n) := \frac{d(n)}{k^{\omega(n)}},
\]
where $d(n)$ is the divisor function and $\omega (n)$ is the number of distinct prime divisors of $n$. The main purpose of this paper is the study of the mean value of the function $D_{k,\omega}(n)$ by using elementary methods.
    }

\keywords{Divisor function, number of distinct prime divisors, mean value.}

\msc{11N37, 11A25, 11N36}

\VOLUME{31}
\YEAR{2023}
\NUMBER{1}
\firstpage{13}
\DOI{https://doi.org/10.46298/cm.10104}

\begin{paper}

\section{Introduction}
Let $k\geq 2$ be a~fixed integer. We recall that $d(n) := \sum_{d\mid n}1$ is the number of divisors of $n$, and $\omega (
n) :=\sum_{p\mid n}1$ is the number of distinct prime divisors
of $n$. We define the function $D_{k,\omega}(n) $ by
\begin{equation}
D_{k,\omega}(n) :=\frac{d(n)}{k^{\omega
(n)}}.
\end{equation}
Notice that for every fixed integer $k\geq 2$, the function $D_{k,\omega}(n) $ is multiplicative and for every prime number $p$ and
every integer $m$ the relation
\begin{equation}
D_{k,\omega}(p^{m}) = \frac{m+1}{k}, \label{11}
\end{equation}
holds. By using \eqref{11}, we get
$$
D_{k,\omega}(n) = \prod\limits_{p^{m}\Vert n}\frac{m+1}{k}
$$
where $p^{m}\Vert n$ means $p^{m}\mid $ $n$ and $p^{m+1}\nmid n$. In the
particular case $k = 2$, the function $D_{2,\omega}(n) $ is
exactly $D(n) = \dfrac{d(n)}{d^{\ast}(
n)},(\text{see~\cite{D.K}})$. For $k\geq 3$, we can
easily check that
\begin{equation}
\sum_{n\leq x}D_{k},_{\omega}(n) \ll _{k}x(\log x)
^{2/k-1}.
\end{equation}
Indeed, for any integer $n$, we have $D_{k},_{\omega}(n) \leq
d(n) \ll _{\varepsilon}n^{\varepsilon}$. Furthermore, the
hypotheses of Shiu's theorem are satisfied; see Theorem $1$ in~\cite{Shui}
and~\cite[p.1]{Nar-Ten}. One gets
$$
\sum_{n\leq x}D_{k},_{\omega}(n) \ll _{k}\frac{x}{\log x}\exp
\Bigl(\sum_{p\leq x}\frac{2}{kp}\Bigr).
$$
Now, by using Lemma $4.63$ in \cite{O.BORD}, it follows that
$$
\sum_{n\leq x}D_{k},_{\omega}(n) \ll _{k}\frac{x}{\log x}\exp
\Bigl(\frac{2}{k}\log (2e^{\gamma}\log x)\Bigr) \ll
_{k}x(\log x) ^{2/k-1}.
$$

\section{Main result}
In this section, we establish two results concerning the mean value of the
function $D_{k,\omega}(n)$. We begin by giving a~weaker result.

\begin{theorem}
Let $k\geq 2$ be a~fixed integer. For all $x\geq 1$ large enough, we have
\[
\sum_{n\leq x}D_{k,\omega}(n) = \frac{x(\log x)^{2/k-1}}{\Gamma
(2/k)}\prod\limits_{p}\Bigl(1-\frac{1}{p}\Bigr) ^{2/k}
\Bigl(1+\frac{2p-1}{kp(p-1) ^{2}}\Bigr) +O(x(\log
x)^{-1}) (\log \log x)^{4/k}.
\]
\end{theorem}

The proof of this result is based on Tulyaganov's theorem; this theorem is
summarized as follows:

\begin{theorem}\label{thm2}
Let $f$ be a~complex valued multiplicative function. Suppose there exists $ z\in\mathbb{C}$, independent of $p$, with $\vert z\vert \leq c_{1}$ and
\begin{enumerate}
\item[a)] $$\sum_{p\leq x}f(p)\log p = zx + O(xe^{-c_{2}\sqrt{\log x}}) $$
\item[b)] $$\sum_{p\leq x}\vert f(p)\vert \log p\ll x
$$
\item[c)] 
$$
\sum_{p\leq x}\sum_{\alpha = 2}^{\infty}\frac{\vert
f(p^{\alpha})\vert \log p^{\alpha}}{p^{\alpha}}\ll (\log \log
x) ^{2}$$
\item[d)] $$\sum_{p}\frac{\vert f(p)\vert ^{2}\log p}{p^{2}}<c_{3}
$$
\end{enumerate}
for some real numbers $c_{1}$, $c_{2}$ and $c_{3}$. Then, for all $x\geq 1$
sufficiently large, we have
\begin{align*}
\sum_{n\leq x}f(n) =
{ }&{ }\frac{x(\log x)^{z -1}} {\Gamma (z )}\prod\limits_{p}\Bigl(1-\frac{1}{p}\Bigr)^{z }\Bigl(1+\sum_{\alpha = 1}^{\infty}\frac{f(p^{\alpha})}{p^{\alpha}}\Bigr) \Bigl\{1+O(\frac{(\log \log x) ^{2}}{\log x}) \Bigr\} \\
&{ }+ O(x(\log x)^{\max (0,\mathrm{Re}\,z-1) -1}) (\log \log x)^{2(A-\max (0,\mathrm{Re}\,z -1))},
\end{align*}
where $A>0$ satisfies
$$
\sum_{u<p\leq v}\vert f(p)\vert p^{-1}\leq A\log (\log
v/\log u) +O(1).
$$
\end{theorem}

\begin{proof}
This theorem is a~consequence of Theorem $4$ in~\cite{Tulyag}, where we take
$g = f$.
\end{proof}

\noindent To complete the demonstration of the main result we have the
following lemmas.

\begin{lemma}\label{lem1}
For any fixed integer $k\geq 2$, we have the estimate
$$
\sum_{p\leq x}\vert D_{k,\omega}(p) \vert \log p\ll
x.
$$
\end{lemma}

\begin{proof}
By Chebyshev's estimates~\cite{Disar}, we have
\[
\sum_{p\leq x}\vert D_{k,\omega}(p) \vert \log p = \frac{2}{k}\sum_{p\leq x}\log p<\frac{2}{k}(1.000081x) \ll x.
\qedhere
\]
\end{proof}

\begin{lemma}\label{lem2}
For any fixed integer $k\geq 2$, there is a~constant $c>0$, such that
$$
\sum_{p\leq x}D_{k,\omega}(p) \log p = \frac{{\footnotesize 2}}{k}x+O(xe^{-c\sqrt{\log x}}).
$$
\end{lemma}

\begin{proof}
We have
$$
\sum_{p\leq x}D_{k,\omega}(p) \log p = \frac{{\footnotesize 2}}{k}\sum_{p\leq x}\log p = \frac{{\footnotesize 2}}{k}\theta (x),
$$
and by Theorem $6.9$ in~\cite{H.L.M and R.C.V}, there is a~constant $c>0$
such that
$$
\theta (x) = x+O(xe^{-c\sqrt{\log x}}),
$$
which implies the desired result.
\end{proof}

\begin{lemma}\label{lem3}
For any fixed integer $k\geq 2$, we have
$$
\sum_{p}\frac{\vert D_{k,\omega}(p) \vert ^{2}}{p^{2}%
}\log p<\infty .
$$
\end{lemma}

\begin{proof}
We first check the inequality $\sum\limits_{m = 2}^{\infty}\dfrac{\log m}{m(m-1)}\leq \log 4$, and using the following
$$
\sum_{p}\frac{\log p}{p^{2}}<\sum_{m = 2}^{\infty}\frac{\log m}{m^{2}}\leq
\sum_{m = 2}^{\infty}\frac{\log m}{m(m-1)},
$$
then we have
\[
\sum_{p}\frac{\vert D_{k,\omega}(p) \vert ^{2}}{p^{2}%
}\log p = \frac{4}{k^{2}}\sum_{p}\frac{\log p}{p^{2}} 
<\frac{4\log 4}{k^{2}}.
\qedhere
\]
\end{proof}

\begin{lemma}\label{lem4}
For any fixed integer $k\geq 2$, we have
$$
\sum_{p\leq x}\sum_{\alpha = 2}^{\infty}\frac{\vert D_{k,\omega}(
p^{\alpha}) \vert \log (p^{\alpha})}{p^{\alpha}}\leq \frac{28}{k}.
$$
\end{lemma}

\begin{proof}
For every integer $k\geq 3$, we write
\begin{align*}
\sum_{p\leq x}\sum_{\alpha = 2}^{\infty}\frac{\vert D_{k,\omega}(
p^{\alpha}) \vert \log (p^{\alpha})}{p^{\alpha}}
&= \frac{1}{k}\sum_{p\leq x}\log p\sum_{\alpha = 2}^{\infty}\frac{\alpha
(\alpha +1)}{p^{\alpha}} \\
&= \frac{1}{k}\sum_{p\leq x}\frac{\log p}{p}\sum_{\alpha = 2}^{\infty}\frac{%
\alpha (\alpha +1)}{p^{\alpha -1}},
\end{align*}
and the infinite series $\sum\limits_{\alpha = 2}^{\infty}\dfrac{\alpha
(\alpha +1)}{p^{\alpha -1}}$ converges to $\dfrac{2}{(
1-1/p) ^{3}}-2$, since
\begin{align*}
\sum_{p\leq x}\sum_{\alpha = 2}^{\infty}\frac{\vert D_{k,\omega}(
p^{\alpha}) \vert \log (p^{\alpha})}{p^{\alpha}}
&= \frac{2}{k}\sum_{p\leq x}\frac{3p^{2}-3p+1}{p(p-1)^{3}}\log p\\
&\leq \frac{28}{k}\sum_{p\leq x}\frac{\log p}{p^{2}}.
\end{align*}
By Lemma $70.1$ in~\cite{HALL -TENENBAUM}, we have $\sum\limits_{p}\dfrac{\log p}{p^{\alpha}}<\dfrac{1}{\alpha -1}$ for all $\alpha >1$, consequently
$$
\sum_{p\leq x}\sum_{\alpha = 2}^{\infty}\frac{\vert D_{k,\omega}(
p^{\alpha}) \vert \log (p^{\alpha})}{p^{\alpha}}<\frac{28}{k}.
$$
Finally, by Lemma \ref{lem1}, \ref{lem2}, \ref{lem3} and \ref{lem4} we have shown that the function $D_{k,\omega}(n) $ satisfies the conditions of Theorem \ref{thm2}. As
we have
$$
\sum_{u<p\leq v}\frac{\vert D_{k,\omega}(p) \vert}{p%
} = \frac{2}{k}\sum_{u<p\leq v}\frac{1}{p}\leq \frac{2}{k}\log \frac{\log v}{%
\log u}+O(1),
$$
then the constant $A$ in Theorem $2$ is $\frac{2}{k}$.
\end{proof}

The next result is improved over the previous one.

\begin{theorem}
Let $k\geq 2$ be a~fixed integer. For all $x\geq 1$ large enough, we have
$$
\sum_{n\leq x}D_{k,\omega}(n) = \frac{x(\log x)^{2/k-1}}{\Gamma
(2/k)}\prod\limits_{p}\Bigl(1-\frac{1}{p}\Bigr) ^{2/k}
\Bigl(1+\frac{2p-1}{k(p-1) ^{2}}\Bigr) +O_{k}(x(\log
x)^{2/k-2}).
$$
\end{theorem}

The demonstration is based on the following lemmas:

\begin{lemma}\label{lem5}
Let $k\geq 2$ be a~fixed integer. For every $s:= \sigma +it\in
\mathbb{C}
$ such that $\sigma >1$ and $L(s,D_{k,\omega}(n))
:= \sum\limits_{n = 1}^{\infty}\dfrac{D_{k,\omega}(n)}{n^{s}}$,
we have
$$
L(s,D_{k,\omega}(n)) = \zeta (s)
^{2/k}L(s,g_{k}),
$$
or $L(s,g_{k}) $ is a~series of Dirichlet absolutely convergent
in the half-plane $\sigma >\frac{1}{2}$.
\end{lemma}

\begin{proof}
If $\sigma >1$, then
\begin{align*}
L(s,D_{k,\omega}(n)) &= \prod\limits_{p}\Bigl(
1+\sum\limits_{\alpha = 1}^{\infty}\dfrac{D_{k,\omega}(p^{\alpha
})}{p^{\alpha s}}\Bigr) \\
&= \prod\limits_{p}\Bigl(1+\sum\limits_{\alpha = 1}^{\infty}\dfrac{\alpha
+1}{kp^{\alpha s}}\Bigr) \\
&= \prod\limits_{p}\Bigl(1+\dfrac{2p^{s}-1}{k(p^{s}-1) ^{2}}\Bigr),
\end{align*}
on the other hand we have
$$
\Bigl(1+\frac{2p^{s}-1}{k(p^{s}-1) ^{2}}\Bigr) = ((
1-p^{-s}) ^{-2/k}) \Bigl(1+\frac{h(s)}{k(
p^{s}-1) ^{2}}\Bigr),
$$
such that
$$
h(s) = (1-p^{-s}) ^{2/k}(kp^{2s}-2(
k-1) p^{s}+k-1) -k(p^{s}-1) ^{2}.
$$
Since
$$
(1-p^{-s}) ^{2/k} = 1-\frac{2}{kp^{s}}-\frac{k-2}{k^{2}p^{2s}}-O\Bigl(\frac{k}{p^{3\sigma}}\Bigr),
$$
he comes
\begin{align*}
h(s) &= \Bigl(1-\frac{2}{kp^{s}}-\frac{k-2}{k^{2}p^{2s}}-O\bigl(\frac{k}{p^{3\sigma}}\bigr)\Bigr) (kp^{2s}-2(
k-1) p^{s}+k-1) -k(p^{s}-1) ^{2} \\
&= 2\bigl(1-\frac{1}{k}\bigr) +O(p^{-\sigma}),
\end{align*}
which implies the announced result.
\end{proof}

\begin{lemma}[\cite{Selb}]\label{lem6}
Let $A>0$. Uniformly for $x\geq 2$ and $z\in
\mathbb{C}
$ such that $\vert z\vert \leq A$, we have
$$
\sum_{n\leq x}\tau _{z}(n) = \frac{x(\log x) ^{z-1}}{%
\Gamma (z)}+O_{A}(x(\log x) ^{\mathrm{Re}\,z-2}).
$$
$\tau _{z}(n) $ is the multiplicative function defined by $\tau
_{z}(p^{\alpha}) = 
\binom{z+\alpha -1}{\alpha}$.
\end{lemma}

\begin{proof}[Proof of Theorem 3]
According to the Lemma \ref{lem5}, we have $D_{k,\omega} = \tau _{2/k}\ast g_{k}$.
Then, by Lemma \ref{lem6}
\begin{align*}
\sum_{n\leq x}D_{k,\omega}(n) &= \sum_{d\leq x}g_{k}(
d) \sum_{m\leq \frac{x}{d}}\tau _{2/k}(m) \\
&= \sum_{d\leq x}g_{k}(d) \Bigl(\frac{x(\log \frac{x}{d}%
) ^{2/k-1}}{d\Gamma (2/k)}+O_{k}\bigl(\frac{x}{d}(
\log \frac{x}{d}) ^{2/k-2}\bigr) \Bigr) \\
&= \sum_{d\leq x}g_{k}(d) \Bigl(\frac{x(\log x)
^{2/k-1}}{d\Gamma (2/k)}+O_{k}((\log x)
^{2/k-2}\log d) \\
&\qquad+O_{k}\bigl(\frac{x}{d}(\log \frac{x}{d})
^{2/k-2}\bigr) \Bigr) \\
&= \frac{x(\log x) ^{2/k-1}}{\Gamma (2/k)}\sum_{d\leq x}\frac{g_{k}(d)}{d}+O_{k}\Bigl(x(\log
x) ^{2/k-2}\sum_{d\leq x}\frac{\vert g_{k}(d)
\vert (1+\log d)}{d}\Bigr).
\end{align*}
The series $L(s,g_{k}) $ is absolutely convergent on the
half-plane $\sigma >\frac{1}{2}$, then for all $\varepsilon >0$
$$
\sum_{d\leq x}\vert g_{k}(d) \vert \ll
_{k,\varepsilon}x^{1/2+\varepsilon},
$$
hence by partial summation
$$
\sum_{d\leq x}\frac{\vert g_{k}(d) \vert (
1+\log d)}{d}\ll _{k,\varepsilon}x^{-1/2+\varepsilon}
$$
and therefore
$$
\sum_{n\leq x}D_{k,\omega}(n) = L(1,g_{k}) \frac{%
x(\log x)^{2/k-1}}{\Gamma (2/k)}+O_{k}(x(\log
x)^{2/k-2}) +O_{k,\omega}(x^{1/2+\varepsilon}).
$$
Which completes the demonstration.
\end{proof}

\section*{Acknowledgments}
The author would like to sincerely thank Professor Olivier Bordellès for his help and interest in this work and Professor Karl Dilcher for his generosity in reviewing this paper.

\EditInfo{29 September, 2019}{28 January, 2020}{Karl Dilcher}

\end{paper}